\documentclass[10.5pt,a4paper]{article}
\usepackage{amsmath,amssymb,mathrsfs,framed,esint,slashed,braket}
\usepackage[colorlinks]{hyperref}
\usepackage{color,multicol,wrapfig}
\usepackage[all]{xy}
\usepackage{graphicx}

\newcommand{\rem}[1]{}

\newcommand{\de}{{\rm d}}

\newcommand{\bp}{{\boldsymbol{p}}}

\newcommand{\bz}{{\mathbf{z}}}

\newcommand{\beq}{\begin{equation}}
\newcommand{\eeq}{\end{equation}}

\newcommand{\ben}{\begin{eqnarray}}
\newcommand{\een}{\end{eqnarray}}

\begin{document}
 
\title{\rm 
Symplectic geometry and Koopman dynamics\\at the quantum-classical interface
}
\sf
\author{CESARE TRONCI}
\date{}

\maketitle

\begin{abstract}\normalsize
Going back to the early days in the history of quantum mechanics, the interaction of quantum and classical systems stands among the most intriguing open questions in science and makes its appearance in  several fields, from physics to chemistry.  Recently, a new perspective on this problem was unfolded by an unprecedented combination of symplectic geometry and Koopman's  formulation of classical mechanics.
\end{abstract}

\bigskip

\begin{multicols}{2}

How does a classical system interact with a quantum system? This simple question is among the most problematic  in science and has been puzzling the community ever since the early discussions among the founding fathers of quantum theory.  Nevertheless, attempts to construct dynamical models of mixed quantum-classical (QC) systems are still ubiquitous in   several fields, from theoretical chemistry to solid state physics. Indeed, the complexity of fully quantum many-body simulations stimulates the search for approximate  models in which part of a quantum system is treated classically while the remainder remains  quantum. The multiscale nature of this problem requires powerful mathematical structures beyond established  methods in semiclassical analysis. 

Despite several efforts, current hybrid QC models suffer from various consistency issues and no general consensus has been reached.
A well-known hybrid QC model is given by the system
\beq\label{mfeqns}
\dot{q}= \partial_p\langle\widehat{H}\rangle
,\qquad 
\dot{p}= -\partial_q\langle\widehat{H}\rangle
,\qquad 
i\hbar\partial_t\psi=\widehat{H}\psi\,, 
\eeq
where $\langle\hat{A}\rangle:=\langle\psi|\hat{A}\psi\rangle$ denotes the expectation value, $\langle\psi_1|\psi_2\rangle$ is the inner product, and $\widehat{H}={H}(q,p,\hat{\sf x},\hat{\sf p})$ is a quantum Hamiltonian operator depending on the classical phase-space coordinates $(q,p)$. Notice that $(\hat{\sf x},\hat{\sf p})$ denote the quantum position and momentum operators such that $[\hat{\sf x},\hat{\sf p}]=i\hbar$. Despite their mathematical appeal, equations \eqref{mfeqns} often fail to produce realistic  results. 

Over the years, it has been recognized that QC coupling requires a probabilistic description in both quantum and classical sectors in such a way that statistical correlation effects are retained in the treatment. This means that, no matter the underlying construction, classical dynamics must ultimately be given in terms of a probability distribution $\rho_c(q,p)$ and quantum evolution in terms of a von-Neumann density operator $\hat\rho$ so that, for example, $\langle\hat{A}\rangle=\operatorname{Tr}(\hat\rho\hat{A})$.  However, blending these essentially different descriptions in such a way to retain statistical QC correlations is far from easy. A new strategy that is recently making its way in   QC coupling  comes from results in pure mathematics, upon blending van Hove's prequantum geometry with Koopman's unitary flows on phase-space.

In 1931  Koopman made the remark that classical mechanics can be formulated as a unitary flow on the Hilbert space of square-integrable complex functions (i.e. wavefunctions, WFs) on phase-space \cite{Koopman}. This  intuitive result  follows by observing that, if $\{\ ,\, \}$ denotes the standard Poisson bracket, one has
\begin{equation}\label{KvN-eq}
i\hbar\partial_{t\!} \chi = i\hbar\{H,\chi\} 
\  \ \implies\ \ 
\partial_t |\chi|^2 = \{H,|\chi|^2\}.
\end{equation}
Then, since the Liouvillian operator $\widehat{L}_H:=i\hbar\{H,\ \} $ is Hermitian, the classical Liouville equation for $\rho_c=|\chi|^2$ can be realized in terms of a unitary evolution on the  Hilbert space $L^2(T^*Q)$. As customary, here we identify  the phase-space with the cotangent bundle $T^*Q$  of the classical configuration manifold $Q$.

The above remark has made many appearances in  the  literature, with several eminent scholars apparently being unaware of the early works by Koopman and von Neumann. These days, the \emph{Koopman-von Neumann equation} (KvN) \eqref{KvN-eq} for the  WF $\chi\!\in\!L^2(T^*Q)$ is attracting increasing attention. 
The idea of using KvN to formulate hybrid QC models goes back to George Sudarshan's work from 1976: if classical mechanics has a Hilbert-space formulation, then one can also take the tensor product of classical and quantum Hilbert spaces, and construct the dynamics of such hybrid QC WFs. However, so far this direction has failed to produce consistent models. Indeed, Koopman's approach involves several challenges  that can hardly be tackled without resorting to powerful mathematics. 

Within an international collaboration \cite{BoGBTr19}, we recently showed how unprecedented insights into hybrid QC dynamics may be uncovered by exploiting the symplectic geometry of Koopman WFs \cite{GBTr20}. After restoring the information on classical phases, a first Koopman-based model was obtained via a partial quantization process upon starting with two fully classical systems. Then, the quantum and classical densities are recovered as momentum maps associated to specific unitary representations corresponding to the classical and quantum motions.

\section*{\normalsize\sf\bfseries Koopman-van Hove prequantum geometry}

As the standard KvN construction failed to produce consistent QC models, we proposed an  alternative approach to Koopman WFs based on \emph{prequantum theory}. First formulated in van Hove's thesis and further developed by Kostant and Souriau, prequantization is a niche branch of symplectic geometry, whose potential in other fields has not been previously realized. Here, we will simply point out that the prequantum evolution differs from the first (KvN) equation in \eqref{KvN-eq} by the insertion of a phase term retaining the information on the Lagrangian function ${L}={\bp\cdot\partial_{\bp} H}-H$. Then, upon restricting to a one-dimensional configuration space for convenience, the KvN equation  becomes 
\begin{equation}\label{KvH-eq}
i\hbar\partial_{t\!} \chi = \widehat{{\cal L}}_H\chi
,\quad\text{with}\quad
\widehat{{\cal L}}_H:= i\hbar\{H,\,\} -(p\partial_p H-H).
\end{equation}
The \emph{prequantum operator} $\widehat{{\cal L}}_H$ has been introduced for later purpose.
At this point,  the polar form $\chi=\sqrt{D}e^{iS/\hbar}$ returns $\de D/\de t=0$ and $\de S/\de t={L}$ along $(\dot{q},\dot{p})=(\partial_pH,-\partial_qH)$. The second equation is a phase-space version of the general Hamilton-Jacobi equation  for Hamilton's principal function. Indeed, unlike the original KvN construction, prequantum theory carries the information on the  classical phase.

Equation \eqref{KvH-eq} was dubbed \emph{Koopman-van Hove} (KvH) equation in \cite{BoGBTr19}. Much insight on \eqref{KvH-eq} is obtained by looking at its canonical Hamiltonian structure, which appears from the variational principle
\beq
\delta\int^{t_2}_{t_1}\!\operatorname{Re}\!\int\!\bar\chi\Big(i\hbar\partial_t\chi-i\hbar\{H,\chi\}+\chi{L}\Big)\de q\de p\,\de t=0
\,.
\label{VP-KvH}
\eeq
Here, the Hamiltonian functional of the system is identified with  the last two integral terms. If we insist that this Hamiltonian must coincide with the physical energy $\int \!\rho_c H\,\de q\de p$, then an integration by parts yields the following expression of the Liouville density in terms of the KvH WF:
\beq
\rho_c=|\chi|^2+\partial_p(p|\chi|^2)+\hbar\operatorname{Im}\{\bar\chi,\chi\}
\,.
\label{momap}
\eeq
While the first term coincides with the original KvN prescription, the remaining two terms appear mysterious. The major breakthrough in \cite{BoGBTr19,GBTr20} was to recognize that the expression above identifies a \emph{momentum map} structure, thereby ensuring that  \eqref{momap} indeed satisfies the classical Liouville equation. We will now give a quick review of this result.
 
\begin{framed}\vspace{-.3cm}
\section*{\normalsize\sf\bfseries What is a momentum map?}
Momentum  maps  are  a  crucial  ingredient in  symplectic  geometry.  Specifically,  the canonical (left) action of a Lie group $G$ with Lie algebra $\mathfrak{g}$ on a symplectic manifold ${\cal S}$ induces a momentum map $J:{\cal S}\to\mathfrak{g}^*$ generalizing Noether’s conserved quantity. Indeed, the latter occurs in the particular case of a symmetry group. When ${\cal S}$ is a linear symplectic space $(V,\Omega)$, one has
\beq\label{momapdef}
\Omega(\xi_V(x),x)=2\langle  J(x),\xi\rangle\,,\quad \forall \xi\in \mathfrak{g}\,,\ \forall x\in V,
\eeq
where $\xi_V$ denotes the infinitesimal $\mathfrak{g}-$action on $V$ and $\langle \,,\rangle$ denotes the duality pairing. Any complex Hilbert space with inner product $\langle \,|\,\rangle$ has the canonical symplectic form  $\Omega(\chi_1,\chi_2)=2\hbar\operatorname{Im}\langle\chi_1|\chi_2\rangle$.
When a Hamiltonian  can be entirely written in terms of a momentum map, this Hamiltonian is called ‘collective' and the dynamics comprises a Lie-Poisson system \cite{GuSt,HoScSt09}. 
Then, the KvH Hamiltonian functional is collective  for the momentum map given by \eqref{momap}. 

\end{framed}

\section*{\normalsize\sf\bfseries Hybrid QC dynamics I: quantization}
But how do we use KvH for devising a mixed quantum-classical model? An important hint is given by the fact that, at least for Hamiltonians of the type $H=T+V$,
 the standard Schr\"odinger equation can be obtained from \eqref{KvH-eq} by simply applying canonical quantization,
 i.e. by 
 enforcing $\partial_p\chi=0$ and replacing $(q,p)\to(\hat{\sf x},-i\hbar\partial_{\sf x})$.
Then, in \cite{BoGBTr19}, a first QC model was obtained by simply starting with a two-particle KvH equation and then quantizing one of them. As a result, the hybrid QC wavefunction $\Upsilon(q,p,{\sf x})$ obeys the quantum-classical wave equation (QCWE)
\beq\label{QCWE}
i\hbar\partial_t\Upsilon=\{i\hbar \widehat{H},\Upsilon\}-\big(p\partial_p\widehat{H}-\widehat{H}\big)\Upsilon
\,,
\eeq
where the Hamiltonian $\widehat{H}$ is the same as  in  \eqref{mfeqns}.  Once again, the RHS identifies a Hermitian operator on the hybrid Hilbert space $\mathscr{H}=L^2(T^*Q\times M)$, where $M$ is the manifold comprising the quantum coordinates ${\sf x}$. 
Thus, the QCWE has again a canonical Hamiltonian structure that is provided by a variational principle analogous to \eqref{VP-KvH}. 
The quantum density matrix $\hat\rho$ and classical Liouville density $\rho_c$ are found as momentum maps corresponding to the actions of quantum unitary operators and van Hove transformations on $\mathscr{H}$, respectively.  More explicitly, we have
\end{multicols}
\begin{framed}\vspace{-.3cm}
\section*{\normalsize\sf\bfseries The van Hove representation}
Consider the group $G=\{(\eta,e^{i\varphi})\in \operatorname{Diff}(T^*Q)\,\circledS\,{\cal F}(T^*Q,S^1)\ |\ \eta^*\theta+\de\varphi=\theta\}$, where $\operatorname{Diff}(T^*Q)$ is the group of diffeomorphisms of $T^*Q$ and ${\cal F}(T^*Q,S^1)$ denotes the space of $S^1-$valued functions on $T^*Q$. Also, $\circledS$ denotes a semidirect product and $\theta$ is the symplectic potential $\theta= p \de q$. The  Lie algebra of $G$ is identified with the space  of scalar functions, endowed with the canonical Poisson bracket. At this point, one  shows that the  operator $-i\hbar^{-1}\widehat{{\cal L}}_H$ generates the unitary \emph{van Hove  representation} of $G$ on  $L^2(T^*Q)$, that is
$
\chi_0(\bz_0)\mapsto \chi_0(\bz_0)e^{{-{i}\varphi(\bz_0)/\hbar}}|_{\bz_0=\eta^{-1}(\bz)}
$.
Then,  the definition \eqref{momapdef} returns the momentum map  \eqref{momap}, whose Lie-Poisson dynamics gives the classical Liouville equation.  
\end{framed}

\begin{multicols}{2}\noindent
\beq\label{quantdens}
\hat\rho({\sf x},{\sf x}')=\int\!\Upsilon(q,p,{\sf x})\bar{\Upsilon}(q,p,{\sf x}')\,\de q\de p
\eeq
\beq\label{classdens}
\rho_c=\int\!\Big(|\Upsilon|^2+\partial_p(p|\Upsilon|^2)+\hbar\operatorname{Im}\{\bar\Upsilon,\Upsilon\}\Big)\de {\sf x}
\eeq
Notice that the quantum density matrix is positive semidefinite. This is the first success	of this approach: as several previous models were unable to retain this property, the current approach represents a substantial step forward. 

As an example, we consider a quantum subsystem given by a $1/2-$spin. In standard Pauli matrix notation, we have $\hat\rho=(\boldsymbol{1}+\mathbf{n}\cdot\widehat{\boldsymbol\sigma})/2$ with $n^2\leq1$. Then, Hamiltonians of the type $\widehat{H}=H_0(q,p)\boldsymbol{1}+H_I(q,p)\widehat{\sigma}_z$ lead to an exactly solvable QCWE \eqref{QCWE} for $\Upsilon\in  L^2(\Bbb{R}^2)\otimes\Bbb{C}^2$, when both $H_0$ and $H_I$ are quadratic. 
\begin{center}
\includegraphics[width=0.49\textwidth]{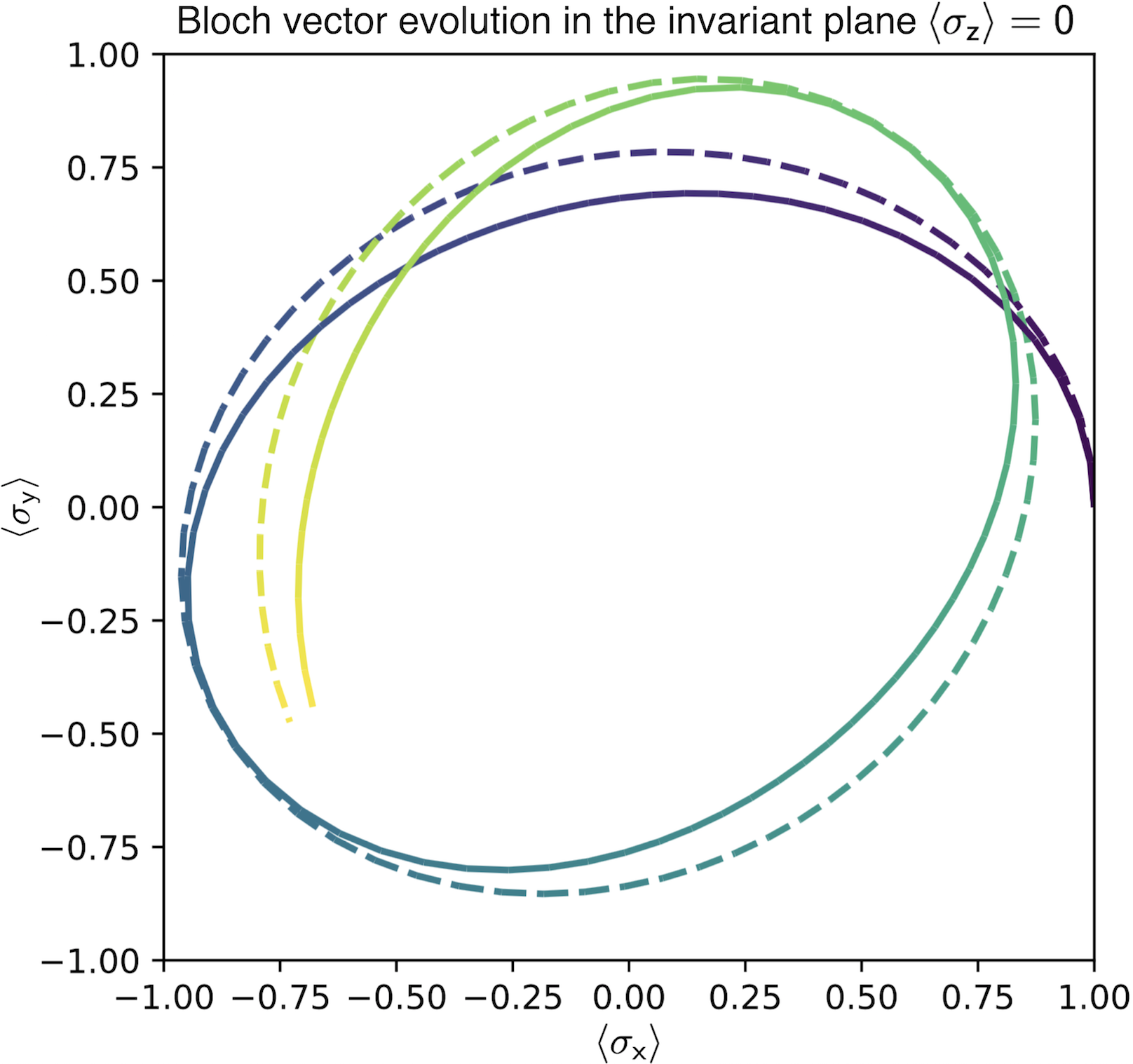}
\end{center}
{
\small\noindent Figure 1. Bloch vector  evolution for    $H_0=(p^2+q^2)/2$ and $H_I=(q^2-p^2)/4+1/2$. The dashed line corresponds to the fully quantum dynamics, while the thick line identifies the QCWE results. The initial time is coloured in purple and the final state (t=10) appears in yellow. The factorized initial condition $\Upsilon_0(q,p)=\chi(q,p)(1,1)/\sqrt{2}$ corresponds to $\mathbf{n}_0=(1,0,0)$ and  $\rho_{c0}=e^{-2H_0}/\pi$.  \emph{Courtesy of Giovanni Manfredi and Icare Morrot-Woisard (IPCMS, Strasbourg)}.
}

\bigskip
\noindent
In this case, the classical density \eqref{classdens} coincides at all times with the results obtained from the fully quantum theory, while a good agreement is also observed for the rotational motion of the  \emph{Bloch vector} $\mathbf{n}=\langle\widehat{\boldsymbol\sigma}\rangle$; see Figure 1. In turn, the QCWE predicts slightly lower levels of  \emph{quantum purity} $\|\hat\rho\|^2=(n^2+1)/2$.

Despite these encouraging results, we  observe that in the case of QC evolution the expression \eqref{classdens} of the classical distribution is generally sign-indefinite. Is the sign of $\rho_c$ preserved by its time evolution? While this is always the case for classical KvH dynamics, a positive answer in the QC context has been found only for an certain infinite family of Hamiltonians $\widehat{H}$  including those discussed above  \cite{GBTr20}. However, no  statement of general validity is currently available. The next section addresses this point by presenting an upgrade model of the QCWE.

\section*{\normalsize\sf\bfseries Hybrid QC dynamics II: $S^1-$symmetry}
So far, things still look rather simple. But how can we restore the positivity of the classical density? A hint is made available by the observation that expressing the KvH equation \eqref{KvH-eq} by using the polar form $\chi=\sqrt{D}e^{iS/\hbar}$ leads to the relation $(\partial_t+\pounds_{X_H})(\de S-\theta)=0$. 
Here, $X_H=(\partial_p H,-\partial_q H)$ is the Hamiltonian vector field and $\pounds$ denotes the Lie derivative. 
Thus, if we could set $\de S=\theta$, the KvH momentum map \eqref{momap} would reduce to the KvN prescription $\rho_c=|\chi|^2$. Unfortunately, things are not that simple because this would introduce challenging topological singularities. In turn, upon using the polar form $\chi=\sqrt{D}e^{iS/\hbar}$, we might want to replace $\de S=\theta$ in the variational principle \eqref{VP-KvH}. This  step has a two-fold effect: a) the Lagrangian becomes trivially $S^1-$invariant, and b) the variational principle \eqref{VP-KvH} returns an alternative formulation of the KvN equation, thereby retaining a positive  density $\rho_c=|\chi|^2$.  Thus,  KvN arises from KvH upon enforcing an $S^1-$symmetry \cite{GBTr22}. 


At this point, one may  apply the same argument to mixed QC dynamics and enforce the $S^1-$symmetry on the QCWE variational principle.   Upon denoting $\widehat{\mathcal{P}}(q,p,{\sf x},{\sf x}')=\Upsilon(q,p,{\sf x})\bar\Upsilon(q,p,{\sf x}')$, one  obtains  
\beq\label{NQCLE}
i\hbar\partial_t \widehat{\mathcal{P}}+i\hbar\operatorname{div} \!\big( \widehat{\mathcal{P}}\big\langle X_{\,\widehat{\mathcal{H}}}\big\rangle\big)=\big[\widehat{\mathcal{H}},\widehat{\mathcal{P}}\big]\,,
\eeq
where  $X_{\widehat{A}}=(\partial_p\widehat{A},-\partial_q\widehat{A}\,)$ and $\langle\widehat{A}\rangle=\operatorname{Tr}(\widehat{A}\widehat{\mathcal{P}})/\rho_c$. Also, $\widehat{\mathcal{H}}= \widehat{H}+\hbar\widehat{F}$, where  $\widehat{F}=\widehat{F}(\widehat{\mathcal{P}},\{\widehat{\mathcal{P}},\widehat{H}\})$ is a prescribed function  \cite{GBTr22}.
Despite its formidable appearance, the nonlinear equation \eqref{NQCLE} appears the first to ensure several consistency properties beyond positivity of the quantum density \eqref{quantdens} and the classical density $\rho_c=\int|\Upsilon|^2\de {\sf x}$. For example,  \eqref{NQCLE} leads to a QC Poincar\'e  invariant extending the classical quantity $\oint_{c(t)}p\de q$. Hence, one  obtains the typical ingredients of symplectic geometry including a symplectic form, a Liouville volume, and Casimir functions.


Starting with the prequantum geometry of Koopman WFs in QC coupling, we have unfolded the emergence of a very rich mathematical structure combining  celebrated concepts in the symplectic geometry of quantum and classical systems. In turn, these geometric structures are capable of accommodating stringent consistency requirements such as quantum and classical positivity, and the reduction to uncoupled quantum and classical dynamics in the absence of a coupling potential.

Computational efforts  are underway, along with the development of fluid closure models. In addition, the present models can be adapted to couple quantum and classical spin systems. An entire new direction is unfolding in front of us and right now we are only scratching the surface of a new ground at the boundary between mathematics, physics, and chemistry.

\section*{\normalsize\sf\bfseries Acknowledgements} 
These results were obtained within an international and interdisciplinary collaboration involving  Denys Bondar (Tulane), Fran\c{c}ois Gay-Balmaz (Paris), Ilon Joseph (Livermore), and Giovanni Manfredi (Strasbourg). The author thanks Paul Bergold, Tom Bridges, and Darryl Holm for their valuable feedback. 
This work was made possible through the support of Grant 62210 from the John Templeton Foundation as well as the Royal Society Grant IES\textbackslash R3\textbackslash203005.

\section*{\sf\bfseries\large Cesare Tronci}

Cesare  is an Associate Professor in mathematics at the University of Surrey with a strong international profile.  His  research interests are in applications of symplectic  geometry, variational fluid models, and multiscale dynamics. Most recently, Cesare became interested in problems of mathematical chemistry.

\end{multicols}

\end{document}